\documentclass[11pt, 4a]{article}

\usepackage{amsmath, amssymb, amsfonts}
\usepackage{authblk}

\title{\bf Epistemic systems and Flagg and Friedman's translation}

\author{Takao Inou\'{e}}

\affil{AI Nutrition Project,\\
Artificial Intelligence Center for Health and Biomedical Research,
National Institutes of Biomedical Innovation, Health and Nutrition (NIBIOHN), Osaka, Japan\footnote{t.inoue@nibiohn.go.jp \\ \indent 
(Personal) takaoapple@gmail.com \\ \indent I prefer my personal mail.}}

\date{October 19, 2023}

\begin{document}
\maketitle

\newtheorem{theorem}{Theorem}[section]
\newtheorem{proposition}{Proposition}[section]
\newtheorem{definition}{Definition}[section]
\newtheorem{corollary}{Corollary}[section]
\newtheorem{lemma}{Lemma}[section]
\newtheorem{convention}{convention}[section]
\newtheorem{remark}{Remark}[section]

\begin{abstract}

In 1986, Flagg and Friedman \cite{ff} gave an elegant alternative proof of the faithfulness of G\"{o}del (or Rasiowa-Sikorski) translation $(\cdot)^\Box$ of Heyting arithmetic $\bf HA$ to Shapiro's epistemic arithmetic $\bf EA$. In \S 2, we shall prove the faithfulness of $(\cdot)^\Box$ without using stability, by introducing another translation from an epistemic system to corresponding intuitionistic system which we shall call \it the modified Rasiowa-Sikorski translation\rm . That is, this introduction of the new translation simplifies the original Flagg and Friedman's proof. In \S 3, we shall give some applications of the modified one for the disjunction property ($\mathsf{DP}$) and the numerical existence property ($\mathsf{NEP}$) of Heyting arithmetic.  In \S 4, we shall show that epistemic Markov's rule $\mathsf{EMR}$ in $\bf EA$ is proved via $\bf HA$. So $\bf EA$ $\vdash \mathsf{EMR}$ and $\bf HA$ $\vdash \mathsf{MR}$ are equivalent. In \S 5, we shall give some relations among the translations treated in the previous sections. In \S 6, we shall give an alternative proof of Glivenko's theorem. In \S 7, we shall propose several 
 (modal-)epistemic versions of Markov's rule for Horsten's modal-epistemic arithmetic $\bf MEA$. And, as in \S 4, we shall study some meta-implications among those versions of Markov's rules in $\bf MEA$ and one in $\bf HA$. Friedman and Sheard gave a modal analogue $\mathsf{FS}$ (i.e. Theorem in \cite{fs}) of Friedman's theorem $\mathsf{F}$ (i.e. Theorem 1 in \cite {friedman}): \it Any recursively enumerable extension of $\bf HA$ which has $\mathsf{DP}$ also has $\mathsf{NPE}$\rm . In \S 8, we shall give a proof of  our \it Fundamental Conjecture \rm $\mathsf{FC}$ which was proposed in Inou\'{e} \cite{ino90a} as follows: $\mathsf{FC}: \enspace \mathsf{FS} \enspace \Longrightarrow \enspace  \mathsf{F}.$ This proof is an example of a new type of proofs. In \S 9, I shall shortly discuss and tell about my philosophy and future plan.
\medskip

\small \it Keywords: \rm Flagg and Friedman's translation, epistemic mathematics, epistemic arithmetic, modal arithmetic, Heyting arithmetic, G\"{o}del translation, .modified Rasiowa-Sikorski translation, Markov's rule, epistemic Markov's rule, disjunction property, numerical existence property,  modal disjunction property, modal numerical existence property, Glivenko's theorem, modal-epistemic arithmetic, Friedman and Sheard's theorem, Friedman's theorem, Fundamental Conjecture for Friedman and Sheard's theorem, equivalence between classical and epistemic theorems, philosophy of mathematics, humanization of mathematics, modality and geometry.

\medskip

\it 2020 Mathematics Subject Classification\rm : 00A30, 03A05, 03B45, 03B70, 03D20, 03D25, 03F03, 03F30, 03F50, 03F55, 03F99, 68Q99, 14A99.


\end{abstract}

\tableofcontents

\section{Introduction, Flagg and Friedman's translation}

In 1986, proposing a translation $(\cdot)_\Gamma^{(E)}$, Flagg and Friedman \cite{ff} gave an elegant alternative proof of the faithfulness of G\"{o}del (or Rasiowa-Sikorski) translation $(\cdot)^\Box$ of Heyting arithmetic $\bf HA$ to Shapiro's epistemic arithmetic $\bf EA$, i.e. an amalgam of Peano arithmeic $\bf PA$ with Lewis's modal logic $\bf S4$. (For details on $\bf EA$, see e.g. Shapiro \cite{shapiro, shapiroE}.) We shall call it \it Flagg and Friedman's translation\rm . As a result, $\bf EA$ is a conservative extension of $\bf HA$. The translation $(\cdot)^\Box$ in \cite{ff} have already enjoyed a rich history established by G\"{o}del \cite{gedel1}, McKinsey and Tarski \cite{mcktar}, Rasiowa and Sikorski \cite{rs}, Maehara \cite{maehara54}, Prawitz and Malmn\"{a}s \cite{pramal}, Shapiro \cite{shapiroE}, Goodman \cite{goodman}, Mints \cite{mints1}, Flagg \cite{flaggE, flagg}, Scedrov \cite{scedrov1, scedrov2, scedrov3, scedrov4}, and others. One may see Troelstra \cite{t33f} for a nice short history until 1984 for $(\cdot)^\Box$. For the further development of it, see Flagg \cite{flaggE, flagg}, Scedrov \cite{scedrov1, scedrov2, scedrov3, scedrov4}. 

In \S 2, we shall prove the faithfulness of $(\cdot)^\Box$ without using stability, by introducing another translation from an epistemic system to corresponding intuitionistic system which we shall call \it the modified Rasiowa-Sikorski translation\rm . That is, this introduction of the new translation simplifies the original Flagg and Friedman's proof. 

In \S 3, we shall give some applications of the modified Rasiowa-Sikorski translation for the disjunction property ($\mathsf{DP}$) and the numerical existence property ($\mathsf{NEP}$) of $\bf HA$. 

In \S 4, we shall show that epistemic Markov's rule $\mathsf{EMR}$ in Shapiro's epistemic arithmetic $\bf EA$, i.e. 
$$\mathsf{EMR}: \enspace \vdash_{\bf EA} \Box \forall x \diamondsuit \exists y R(x, y) \enspace \Rightarrow \enspace \vdash_{\bf EA} \Box \forall x \exists y R(x, y),$$
is proved via Heyting arithmetic $\bf HA$: the deletion of $\diamondsuit$, i.e. $\neg \Box \neg$ from $\Box \forall \diamondsuit \exists y R(x, y)$ can be reduced to Markov's rule $\mathsf{MR}$ in $\bf HA$, i.e.
 $$ \mathsf{MR}: \enspace \vdash_{\bf HA} \forall x \neg \neg \exists y R(x, y) \enspace \Rightarrow \enspace \vdash_{\bf HA} \forall x \exists y R(x, y),$$
where $R(x, y)$ is $f(x, y) = 0$ for some binary primitive recursive function symbol $f$. This result was announced in Inou\'{e} \cite{ino90b}. Flagg and Friedman \cite{ff} showed that $\mathsf{EMR}$ implies $\mathsf{MR}$. Our result for which we give two different proofs in \S 4 would be the first example of the reduction of epistemic to intuitionistic properties. In the first proof of it, we make use of Flagg and Friedman's translation introduced in \cite{ff}, while in the second one we apply to the problem the conservativeness result (Mints \cite{mints1} and Goodman \cite{goodman}) that $\bf EA$ is a conservative extension of $\bf HA$. So $\bf EA$ $\vdash \mathsf{EMR}$ and $\bf HA$ $\vdash \mathsf{MR}$ are equivalent.

In \S 5, we shall give some interesting relations among the translations treated in the 
previous sections.

In \S 6, we shall give an alternative proof of Glivenko's theorem. Glivenko's pioneer paper \cite{glivenko} can be regarded as a precursor of G\"{o}del \cite{gedel1} in which his celebrated reduction of classical to intuitionistic formal systems is carried out. The celebrated Glivenko's theorem will be proved by making use of Flagg and Friedman's translation. This proof-method used in the proof is a generalization of the well-known one: in other words, one instance of Flagg and Friedman's translation contains G\"{o}del translation from classical to intuitionistic propositional (predicate) logics. We can find a model-theoretic proof of Glivenko's theorem as an application of modal logic e.g. in Fitting \cite[p. 552]{fitting}. However, our treatment for it is proof-theoretic.

In \S 7, we shall present some results about Horsten's modal-epistemic arithmetic 
and modal-epistemic Markov's rules. In Shapiro's epistemic arithmetic $\bf EA$, the operator $\Box$ (which is Shapiro's $K$, properly speaking) was introduced as provability (or knowability): that is, $\Box A$ intuitively means that $A$ is provable. But we can further analyse $\Box A$, for example as ``it is possible that we have a proof of $A$". So the provability-operator $\Box$ can be divided to two parts, namely the possibility-operator $\Diamond$ and the proof-operator $P$.  In this context, we can regard $\Box A$ as $\Diamond PA$.  This is the motivation of Leon Horsten to introduce his modal-epistemic arithmetic $\bf MEA$ (see Horsten \cite{horsten1}). In the formulation of $\bf MEA$, $\Box$ and $P$ are taken as primitive operators instead of only $\Box$ in the sense of Shapiro. The possibility-operator  is defined as usual. In \S 7, we shall propose several (modal-)epistemic versions of Markov's rule for Horsten's modal-epistemic arithmetic $\bf MEA$. And, as in \S 4, we shall study some meta-implications among those versions of Markov's rules in $\bf MEA$ and one in $\bf HA$.

Friedman and Sheard gave a modal analogue $\mathsf{FS}$ (i.e. Theorem in \cite{fs}) of Friedman's theorem $\mathsf{F}$ (i.e. Theorem 1 in \cite {friedman}): \it Any recursively enumerable extension of Heyting arithmetic $\bf HA$ which has the disjunction property $\mathsf{DP}$ also has the numerical existence property $\mathsf{NPE}$\rm . In \S 8, we shall give a proof of  our \it Fundamental Conjecture \rm $\mathsf{FC}$ which was proposed in Inou\'{e} \cite{ino90a} as follows: $\mathsf{FC}: \enspace \mathsf{FS} \enspace \Longrightarrow \enspace  \mathsf{F}.$  This proof is an example of a new type of proofs.

In \S 9, I shall shortly discuss and tell about my philosophy and future plan.

From here to the end of this section, we shall deal with preliminaries for the later sections.

Throughout this paper, the logical symbols $\vee$ (disjunction) , $\wedge$ (conjunction), $\supset$ (implication), $\exists$ (existential quantifier), $\forall$ (universal quantifier) and the logical constant $\bot$ (falsum) are taken as primitives for intuitionistic first-order calculus with equality $\bf IQC$. We add to these the modal operator $\Box$ for the corresponding epistemic system $\bf EQC$, i.e. modal first-order calculus $\bf S4$ with equality. We write $\neg A$ for $(A \supset \bot)$. A system of natural deduction for $\bf IQC$ is formulated as follows (For natural deduction, see e.g. Prawitz \cite{prawitz} or van Dalen \cite{vdalen}):

$$\mbox{$\wedge$ I} \enspace \frac{\enspace A_0 \enspace \enspace A_1 \enspace }{A_0 \wedge A_1} \enspace \enspace \enspace \enspace \enspace \enspace \mbox{$\wedge$ E} \enspace \frac{\enspace A_0 \wedge A_1 \enspace }{\enspace A_i \enspace}$$

$$\mbox{$\vee$ I} \enspace \frac{\enspace A_i \enspace}{\enspace A_0 \vee A_1 \enspace } \enspace \enspace \enspace \enspace \enspace \enspace \mbox{$\vee$ E} \enspace \frac{\enspace A_0 \vee A_1 \enspace \enspace \begin{array}{c} [A_0] \\ B \end{array} \enspace \enspace \begin{array}{c} [A_1] \\ B \end{array} \enspace }{\enspace B \enspace}$$

$$\mbox{$\supset$ I} \enspace \frac{\enspace \begin{array}{c} [A] \\ B \end{array} \enspace }{\enspace A \supset B \enspace } \enspace \enspace \enspace \enspace \enspace \enspace \mbox{$\supset$ E} \enspace \frac{\enspace A \enspace \enspace A \supset B \enspace }{\enspace B \enspace}$$

$$\bot_{\rm I} \enspace \frac{\enspace \bot \enspace }{\enspace A \enspace }$$
$$\mbox{$\exists$ I} \enspace \frac{\enspace A [t/x] \enspace}{\enspace \exists x A \enspace } \enspace \enspace \enspace \enspace \enspace \enspace \mbox{$\exists$ E} \enspace \frac{\enspace \exists x A  \enspace \enspace \begin{array}{c} [A] \\ B \end{array} \enspace }{\enspace B \enspace}$$

$$\mbox{$\forall$ I} \enspace \frac{\enspace A \enspace}{\enspace \forall x A \enspace } \enspace \enspace \enspace \enspace \enspace \enspace \mbox{$\forall$ E} \enspace \frac{\enspace \forall x A  \enspace }{\enspace A [t/x] \enspace}$$

\medskip

\bf Remark\rm . The usual restriction on the above is adopted.

The rule for equality are the following.

$$\mbox{$\bf Eq$} \enspace t = t \enspace \enspace \enspace \enspace \enspace \enspace \mbox{$\bf Sub$} \enspace \frac{\enspace s = t  \enspace \enspace \enspace A(s) \enspace }{\enspace A(t) \enspace}$$

The rules for $\bf EQC$ is obtained by dropping the rule $\bot_{\rm I}$ for the above set of rules and adding the following.

$$\mbox{$\Box$ I} \enspace \frac{\enspace A \enspace }{\enspace \Box A \enspace}$$

\medskip

\noindent (provided all open assumptions are of the form $\Box B$)

$$\mbox{$\Box$ E} \enspace \enspace \frac{\enspace  \Box A  \enspace }{\enspace A \enspace} \enspace \enspace \enspace \enspace \enspace \enspace \mbox{$\bot_{\rm C}$} \enspace \frac{\enspace \begin{array}{c} [\neg A] \\ \bot \end{array} \enspace }{\enspace A \enspace}$$

\bf Remark\rm . As far as I am concerned, Gentzen-style formulations for the above, e.g. Kleene's $\bf G3$ for $\bf IQC$, Ohnishi and Matsumoto's calculus for $\bf EQC$ practically are more convenient for checking derivations. See Kleene \cite{kleene1} and  Ohnishi and Matsumoto \cite{om}. Also see Chagrov and Zakharyaschev \cite{cz}, Hughes and Cresswell \cite{hc1, hc2}, Sch\"{u}tte \cite{schutteV, schutte2} and Takeuti \cite{takeuti}. We also do not forget Ketonen \cite{ketonen} (cf. Negri and von Plato \cite{NegriPlato2022}  and von Boguslawsk \cite{Boguslawski2010}). For a contradiction calculus, Inou\'{e} 
\cite{Inoue1989} gave a Ketonen-type Gentzen system.  In a direction of computer science, refer to Inou\'{e} \cite{Inoue2003, Inoue2023miz},  Inou\'{e} and Hanaoka  \cite{InoueHanaoka2022} and Inou\'{e} and Miwa \cite{InoueMiwa2023} as studies for Mizar system and 
Mizar Mathematical Library MML. 

If $S$ is one of the formal systems and $A_1, \dots, A_n, B$ are formulas of $S$, we write 
$$A_1, \dots, A_n \vdash_{S} B$$
to indicate that there is a derivations in $S$ of $B$ with open assumptions among $A_1, \dots, A_n$. We sometimes write $S \vdash A$ for $\vdash_{S} A$. Also we write $A \dashv \vdash_{S} B$ if $A \vdash_{S} B$ and $B \vdash_{S} A$.

By $\bf EA$ we denote the so-called Shapiro's epistemic arithmetic introduced by Shapiro \cite{shapiro} and independently Reinhardt \cite{reinhardt}. The constants of $\bf EA$ are the usual $0$, $s$, $+$, $\cdot \enspace$. The axioms and rules of $\bf EA$ are those of $\bf EQC$ together with the usual axioms for zero $0$ and successor $s$, the defining equations for $+$, $\cdot$ and the scheme of induction for all formulas of the language. In other words, $\bf EA$ is Peano arithmetic $\bf PA$ $+$ $\bf S4$. For Heyting arithmetic $\bf HA$, we take the usual formulation corresponding to $\bf EA$ (see e.g. Troelstra and van Dalen \cite{tvd}).

Here we recall the definition of G\"{o}del (or Rasiowa-Sikorski) translation $(\cdot)^\Box$ from $\bf IQC$ $(\bf HA)$ to $\bf EQC$ $(\bf EA)$.

\begin{definition} \rm For each formula $A$ of $\bf IQC$ $(\bf PA)$, we define, by induction on the complexity of $A$, the formula $A^\Box$ of $\bf EQC$ $(\bf EA)$ as follows:
\begin{enumerate}
\renewcommand{\theenumi}{(\roman{enumi})}
\renewcommand{\labelenumi}{\theenumi}
\item  If $A = \bot$ or $s=t$, then $A^\Box = A$. 
\item  If $A$ is an atomic formula different from $\bot$ and $s=t$, then $A^\Box = \Box A$. 
\end{enumerate}
\begin{tabbing}
\hspace*{.05in} \= \hspace{.4in} \= \hspace{.7in} \= \hspace{.3in} \= \kill

\> (iii) \> $(A \# B)^\Box$ \> $=$ \> $A^\Box \# B^\Box$, $\enspace$ $\# \in \{\vee, \wedge\}$.\\

\> (vi) \> $(A \supset B)^\Box$ \> $=$ \> $\Box (A^\Box \supset B^\Box)$. \\

\> (v) \> $(\exists x A)^\Box$ \> $=$ \> $\exists x A^\Box$. \\

\> (vi) \> $(\forall x A)^\Box$ \> $=$ \> $\Box \forall x A^\Box$.

\end{tabbing}

\end{definition} 

\begin{definition} \rm A formula $A$ of a modal formal system $S$ is called \it $S$-stable \rm if 
$$S \vdash A \leftrightarrow \Box A.$$
\end{definition}

\begin{lemma} Let $S$ be $\bf EQC$.

$(1)$ Each formula of the form $\Box A$ is $S$-stable.

$(2)$ If $A$ and $B$ are $S$-stable, then so are $(A \vee B)$ and $(A \wedge B)$.

$(3)$ The following is a permissible inference rule of $S$:

$$\mbox{$\Box$ $\rm I'_{\it S}$} \enspace \enspace \frac{\enspace  A  \enspace }{\enspace \Box A \enspace}$$

\noindent $($provided all open assumptions are $S$-stable$)$.

$(4)$ For any terms $s$ and $t$, the formula $(s = t)$ is $S$-stable.

$(5)$ If $A$ is $S$-stable, then so is $\exists x A$.

\end{lemma}

\begin{lemma} For each formula $A$ of $\bf EQC$, $A^\Box$ is $\bf EQC$-stable.

\end{lemma}

After \cite{ff}, we shall introduce a convenient notation 
$$\neg_{E} A \enspace =_{def.} \enspace A \supset E$$
for any two formulas $A$ and $E$ of intuitionistic formal systems. We write $\neg ^n_E A$ for 
$$\underbrace{\neg_E \cdots \neg_E}_{\mbox{$n$ times}} A.$$

The following lemma is useful.

\begin{lemma} \rm (Flagg and Friedman \cite{ff}) \it Let $S$ be $\bf IQC$. Let $\Gamma$ be a finite set of formulas of $S$ and let $A$, $B$, $E$ be formulas of $S$. Then we have the following.

$(1)$ \enspace $A \dashv \vdash_S \neg^2_A A$.

$(2)$ \enspace $A \vdash_S \neg^2_E A$.

$(3)$ \enspace $\neg_E A \dashv \vdash_S \neg^3_E A$.

$(4)$ \enspace If $A \in \Gamma$, then 
$$A \dashv \vdash_S \bigwedge_{E \in \Gamma} \neg^2_E A.$$

$(5)$ \enspace $\neg^2_E (A \vee B) \dashv \vdash_S \neg^2_E (\neg^2_E A \vee \neg^2_E B).$

$(6)$ \enspace $\neg^2_E (A \wedge B) \dashv \vdash_S \neg^2_E A \wedge \neg^2_E B$.

$(7)$ \enspace If $A \supset B \in \Gamma$, then 
$$A \supset B \dashv \vdash_S \bigwedge_{E \in \Gamma} (\neg^2_E A \supset \neg^2_E B).$$

$(8)$ \enspace $\neg^2_E \exists x A \dashv \vdash_S \neg^2_E \exists x \neg^2_E A$.

$(9)$ \enspace The following rule is permissible in $S$:
$$\frac{\enspace A_1, \dots, A_n \vdash_{S} B \enspace}{\enspace \neg^2_E A_1, \dots, \neg^2_E A_n \vdash_{S} \neg^2_E B \enspace.}$$

\end{lemma}

Now we shall recall the definition of \it Flagg and Friedman's translation \rm $(\cdot)_\Gamma^{(E)}$.

\begin{definition} \rm (Flagg and Friedman \cite{ff}) Let $\Gamma$ be a finite set of formulas of $\bf IQC$ and let $E \in \Gamma$. For each formula $A_\Gamma^{(E)}$ of $\bf IQC$ as follows.

$(1)$ \enspace $A_\Gamma^{(E)} \equiv \neg_E \neg_E A$ for each atomic formula $A$.

$(2)$ \enspace $(A \vee B)_\Gamma^{(E)} \equiv \neg_E \neg_E (A_\Gamma^{(E)} \vee B_\Gamma^{(E)}).$

$(3)$ \enspace $(A \wedge B)_\Gamma^{(E)} \equiv (A_\Gamma^{(E)} \vee B_\Gamma^{(E)}).$

$(4)$ \enspace $(A \supset B)_\Gamma^{(E)} \equiv (A_\Gamma^{(E)} \supset B_\Gamma^{(E)})$.

$(5)$ \enspace $(\Box A)_\Gamma^{(E)} \equiv \neg_E \neg_E \bigwedge_{C \in \Gamma} A_\Gamma^{(C)}$.

$(6)$ \enspace $(\exists x A)_\Gamma^{(E)} \equiv \neg_E \neg_E \exists x A_\Gamma^{(E)}$.

$(7)$ \enspace $(\forall x A)_\Gamma^{(E)} \equiv \forall x A_\Gamma^{(E)}$.
\end{definition}

We write $A_\Gamma^E$ for $A_\Gamma^{(E)}$ if no confusion arises. We also write $$\bigwedge_{C} A_\Gamma^C$$ for $$\bigwedge_{C \in \Gamma} A_\Gamma^{(C)}$$ if no confusion arises.

The translation $(\cdot)_\Gamma^{(E)}$ has the following favorable properties.

\begin{lemma} \rm (\cite{ff}) \it For any formulas $B$ of $\bf EQC$ and for any finite set $\Gamma$ of formulas of $\bf IQC$, and any $E \in \Gamma$,

$(1)$ \enspace $(\neg B)_\Gamma^E \dashv \vdash_{\bf IQC} \neg_E B^E_\Gamma.$

$(2)$ \enspace $\neg_E^2  B_\Gamma^E \dashv \vdash_{\bf IQC} B^E_\Gamma.$

\end{lemma}

We see the soundness of $(\cdot)_\Gamma^{(E)}$ by induction on derivations.

\begin{lemma} \rm (\cite{ff}) \it Let $A_1, \dots, A_n, B$ be formulas of $\bf EQC$. If  
$$A_1, \dots, A_n \vdash_{\bf EQC} B,$$
then for any finite set $\Gamma$ of formulas of $\bf IQC$, and any $C \in \Gamma$,
$$A_{1 \Gamma}^{(C)}, \dots, A_{n \Gamma}^{(C)} \vdash_{\bf IQC} B_\Gamma^{(C)}.$$
\end{lemma}

\begin{theorem} \rm (\cite{ff}) \it Let $B$ be a formula of $\bf EQC$. If  
$$\vdash_{\bf EQC} B,$$
then for any finite set $\Gamma$ of formulas of $\bf IQC$, and any $C \in \Gamma$,
$$\vdash_{\bf IQC} B_\Gamma^{(C)}.$$
\end{theorem}

\begin{lemma} \rm (\cite{ff}) \it Let $A_1, \dots, A_n, B$ be formulas of $\bf EA$. If  
$$A_1, \dots, A_n \vdash_{\bf EA} B,$$
then for any finite set $\Gamma$ of formulas of $\bf IQC$, and any $C \in \Gamma$,
$$A_{1 \Gamma}^{(C)}, \dots, A_{n \Gamma}^{(C)} \vdash_{\bf HA} B_\Gamma^{(C)}.$$
\end{lemma}

\begin{theorem} \rm (\cite{ff}) \it Let $B$ be a formula of $\bf EA$. If  
$$\vdash_{\bf EA} B,$$
then for any finite set $\Gamma$ of formulas of $\bf HA$, and any $C \in \Gamma$,
$$\vdash_{\bf HA} B_\Gamma^{(C)}.$$
\end{theorem}

\begin{lemma} \rm (\cite{ff}) \it Let $A$ be a formula $\bf IQC$ and suppose that $\Gamma$ consists of all subformulas of $A$. The we have  
$$A \dashv \vdash_{\bf IQC} \bigwedge_{C \in \Gamma} (A^\Box)_\Gamma^{(C)}.$$
\end{lemma} 

\it Flagg and Friedman's uniform method \rm is a series of procedures using Lemma 1.7 as the core of the method to give an elegant proof of the following theorem. To prove Lemma 1.7, Flagg and Friedman made use of stability defined in Definition 1.2.

\begin{theorem} \rm (The soundness and faithfulness of $(\cdot)^\Box$) (cf. \cite{ff}) \it For each formula $A$ of $\bf IQC$ $(\bf HA)$, 
$$\vdash_{\bf IQC (HA)} A \enspace \Longleftrightarrow \enspace \enspace \vdash_{\bf EQC (EA)} A^\Box.$$
\end{theorem} 

The original result of Theorem 1.3 was proved by Rasiowa and Sikorski \cite{rs} and independently by Prawitz and Malmn\"{a}s \cite{pramal} for intuitionistic predicate logic.

\it Proof of Theorem 1.3 by Flagg and Friedman's uniform method\rm . The implication $\Rightarrow$ is proved by induction on derivations. For the converse, suppose $\bf EQC (EA)$ $\vdash A^\Box$. Let $\Gamma$ be the set of all subformulas of $A$. By Theorems 1.1 and 1.2, we have 
$$\vdash_{\bf IQC} \bigwedge_{C \in \Gamma} (A^\Box)_\Gamma^{(C)}.$$
Then apply Lemma 1.7 to get $\vdash_{\bf IQC (HA)} A$. $\Box$

As an immediate consequence of Theorem 1.3, we have:

\begin{corollary} $\bf EA$ $($$\bf EQC$$)$ is a conservative extension of $\bf HA$$($$\bf IQC$$)$.

\end{corollary} 

Flagg and Friedman's uniform method is based on Funayama's Theorems (see Funayama \cite{funayama} and e.g. Johnstone \cite[p. 253]{johnstone}, Flagg \cite{flaggC}): for any complete Heyting algebra cHa, there exists a complete Boolean algebra cBa and embedding cHa $\rightarrow$ cBa which preserves finite meets and arbitrary joins. 

\section{The modified Rasiowa and Sikorski's translation}

In this section, we shall introduce \it the modified Rasiowa-Sikorski translation $($or interpretation$)$ $(\cdot)^{\Box_{\mathbf{RS}}}$ \rm in order to develop Friedman's uniform method without stability. 

\begin{definition} \rm For each formula $A$ of $\bf IQC$ $(\bf PA)$, we define, by induction on the complexity of $A$, the formula $A^{\Box_{\mathbf{RS}}}$ of $\bf EQC$ $(\bf EA)$ as follows:
\begin{enumerate}
\renewcommand{\theenumi}{(\roman{enumi})}
\renewcommand{\labelenumi}{\theenumi}
\item  If $A = \bot$ or $s=t$, then $A^{\Box_{\mathbf{RS}}} = A$. 
\item  If $A$ is an atomic formula different from $\bot$ and $s=t$, then $A^{\Box_{\mathbf{RS}}} = \Box A$. 
\end{enumerate}
\begin{tabbing}
\hspace*{.05in} \= \hspace{.4in} \= \hspace{.9in} \= \hspace{.3in} \= \kill

\> (iii) \> $(A \vee B)^{\Box_{\mathbf{RS}}}$ \> $=$ \> $\Box A^{\Box_{\mathbf{RS}}} \vee \Box B^{\Box_{\mathbf{RS}}}$.\\

\> (iv) \> $(A \wedge B)^{\Box_{\mathbf{RS}}}$ \> $=$ \> $A^{\Box_{\mathbf{RS}}} \wedge B^{\Box_{\mathbf{RS}}}$.\\

\> (v) \> $(A \supset B)^{\Box_{\mathbf{RS}}}$ \> $=$ \> $\Box (\Box A^{\Box_{\mathbf{RS}}} \supset \Box B^{\Box_{\mathbf{RS}}})$. \\

\> (vi) \> $(\exists x A)^{\Box_{\mathbf{RS}}}$ \> $=$ \> $\exists x \Box A^{\Box_{\mathbf{RS}}}$. \\

\> (vii) \> $(\forall x A)^{\Box_{\mathbf{RS}}}$ \> $=$ \> $\forall x A^{\Box_{\mathbf{RS}}}$.

\end{tabbing}

\end{definition} 

A version of $(\cdot)^{\Box_{\mathbf{RS}}}$ has been studied in Prawitz and Malmn\"{a}s \cite{pramal}. Our $(\cdot)^{\Box_{\mathbf{RS}}}$ is more close to the original G\"{o}del's translation for propositional logic (see G\"{o}del \cite{gedel1}) than $(\cdot)^\Box$ in Flagg and Friedman \cite{ff}.

We further note:

\begin{proposition} For each formula $A$ of $\bf IQC$ $(\bf PA)$, we have 
$$\vdash_{\bf EQC (EA)} A^\Box \leftrightarrow A^{\Box_{\mathbf{RS}}}.$$

\end{proposition}

Proof. Easy by induction on the complexity of $A$. $\Box$

Although we do not need the following lemma for proving our main theorem, we show it in comparison with the stability of $(\cdot)^\Box$.

\begin{lemma} $ $ 

$(1)$ For each formula $A$ of $\bf EQC$, $A^{\Box_{\mathbf{RS}}}$ is $\mathbf{EQC}$-stable, if $A$ is not of the form $\forall x B$.

$(2)$ For each formula $A$ of $\bf EQC$, $A^{\Box_{\mathbf{RS}}}$ is $\mathbf{EQC}$-stable, if the Barcan formulas $($i.e. formulas of the form $\forall x \Box B \supset \Box \forall x B$$)$ are added to $\bf EQC$ as axioms.

\end{lemma}

By induction on derivation, we can easily prove the following soundness of $(\cdot)^{\Box_{\mathbf{RS}}}$. 
Or use Proposition 2.1 and Lemma 1.1, simply.

\begin{theorem} \rm (Soundness) \it For each formula $A$ of $\bf IQC$ $(\bf HA)$, 
$$\vdash_{\bf IQC (HA)} A \enspace \Longrightarrow \enspace \enspace \vdash_{\bf EQC (EA)} A^{\Box_{\mathbf{RS}}}.$$
\end{theorem} 

Our main concern in this section is to prove the following faithfulness of $(\cdot)^{\Box_{\mathbf{RS}}}$ using Flagg and Friedman's uniform method \it without stability\rm .

\begin{theorem}  \rm (Faithfulness) \it For each formula $A$ of $\bf IQC$ $(\bf HA)$, 
$$\vdash_{\bf EQC (EA)} A^{\Box_{\mathbf{RS}}} \enspace \Longrightarrow \enspace \enspace \vdash_{\bf IQC (HA)} A.$$
\end{theorem} 

The following lemma is essential for Flagg and Friedman's uniform method with respect to $(\cdot)^{\Box_{\mathbf{RS}}}$.

\begin{lemma} Let $A$ be a formula $\bf IQC$ and suppose that $\Gamma$ consists of all subformulas of $A$. The we have  
$$A \dashv \vdash_{\bf IQC} \bigwedge_{C \in \Gamma} (A^{\Box_{\mathbf{RS}}})_\Gamma^{(C)}.$$
\end{lemma} 

Proof. For brevity, we write $A^\Box$ for $A^{\Box_{\mathbf{RS}}}$ only in this proof. We often write $A^C$ for $A^C_\Gamma$. 

Case 1: ($A \equiv \bot$ or $(s = t)$) Apply Lemma 1.3.(1).

Case 2: ($A$ is an atomic formula $P$ which does not fall in Case 1.) We see
$$P^{\Box C} \equiv \neg^2_C \neg^2_P P \enspace \mbox{and} \enspace P \dashv \vdash_{\bf IQC} \neg^2_P P.$$
Thus we have
$$\neg^2_C P \dashv \vdash_{\bf IQC} P^{\Box C} \enspace \mbox{so} \enspace \bigwedge_{C} \neg^2_C P \dashv \vdash_{\bf IQC} \bigwedge_{C} P^{\Box C}.$$
By Lemma 1.3.(6), we get 
$$P \dashv \vdash_{\bf IQC} \bigwedge_{C} P^{\Box C}.$$

Case 3: ($A \equiv A_0 \vee A_1$) We see 
$$A^{\Box C} \equiv \neg^2_C (\neg^2_C  \bigwedge_{E} A_0^{\Box E} \vee \neg^2_C \bigwedge_{E} A_1^{\Box E}).$$
Using Lemma 1.3.(5), we have 
$$\neg^2_C (\bigwedge_{E} A_0^{\Box E} \vee \bigwedge_{E} A_1^{\Box E}) \dashv \vdash_{\bf IQC} A^{\Box C}.$$
By induction hypothesis (for short, I.H.) and Lemma 1.3.(9), it follows that 
$$\neg^2_C (A_0 \vee A_1) \dashv \vdash_{\bf IQC} A^{\Box C},$$
so
$$\bigwedge_{C} \neg^2_C (A_0 \vee A_1) \dashv \vdash_{\bf IQC} \bigwedge_{C} A^{\Box C}.$$
By Lemma 1.3.(3), we get 
$$A_0 \vee A_1 \dashv \vdash_{\bf IQC} \bigwedge_{C} A^{\Box C}.$$

Case 4: ($A \equiv A_0 \wedge A_1$) Then $A^{\Box C} \equiv A_0^{\Box C} \wedge A_1^{\Box C}$. We see that
$$\bigwedge_{C} (A_0^{\Box C} \wedge A_1^{\Box C}) \dashv \vdash_{\bf IQC} \bigwedge_{C} A^{\Box C},$$
namely
$$\bigwedge_{C} A_0^{\Box C} \wedge \bigwedge_{C} A_1^{\Box C} \dashv \vdash_{\bf IQC} \bigwedge_{C} A^{\Box C},$$
so by I.H., 
$$A_0 \wedge A_1 \dashv \vdash_{\bf IQC} \bigwedge_{C} A^{\Box C}.$$

Case 5: ($A \equiv A_0 \supset A_1$) We first see 
$$A^{\Box C} \equiv \neg^2_C \bigwedge_{E} (\neg^2_E  \bigwedge_{F} A_0^{\Box F} \supset \neg^2_E \bigwedge_{F} A_1^{\Box F}).$$
By I.H., it follows that
$$\neg^2_C \bigwedge_{E} (\neg^2_E A_0 \supset \neg^2_E A_1) \dashv \vdash_{\bf IQC} A^{\Box C}.$$
Applying Lemma 1.3.(7), we obtain
$$\neg^2_C (A_0 \supset A_1) \dashv \vdash_{\bf IQC} A^{\Box C}.$$
Then we have 
$$\bigwedge_{C} \neg^2_C (A_0 \supset A_1) \dashv \vdash_{\bf IQC} \bigwedge_{C} A^{\Box C}.$$
Now we use Lemma 1.3.(4) to get 
$$A_0 \supset A_1 \dashv \vdash_{\bf IQC} \bigwedge_{C} A^{\Box C},$$
since $A_0 \supset A_1 \in \Gamma$.

Case 6: ($A \equiv \exists x A_0$) We observe that 
$$\neg^2_C \exists x \neg^2_C \bigwedge_{E} A_0^{\Box E} \dashv \vdash_{\bf IQC} A^{\Box C}.$$
By Lemma 1.3.(8) and I.H., we have
$$\neg^2_C \bigwedge_{E} A_0^{\Box E} \dashv \vdash_{\bf IQC} A^{\Box C}.$$
Then we have 
$$\bigwedge_{C} \neg^2_C \bigwedge_{E} A_0^{\Box E} \dashv \vdash_{\bf IQC} \bigwedge_{C} A^{\Box C}.$$
Then apply Lemma 1.3.(4) to obtain 
$$\exists x A_0 \dashv \vdash_{\bf IQC} \bigwedge_{C} A^{\Box C},$$
since $A_0 \in \Gamma$.

Case 7: ($A \equiv \forall x A_0$) By the definition, we see that 
$$\forall x A_0^{\Box C} \dashv \vdash_{\bf IQC} A^{\Box C},$$
thus
$$\bigwedge_{C} \forall x A_0^{\Box C} \dashv \vdash_{\bf IQC} \bigwedge_{C} A^{\Box C}.$$
In $\bf IQC$, it immediately follows that 
$$\forall x \bigwedge_{C} A_0^{\Box C} \dashv \vdash_{\bf IQC} \bigwedge_{C} A^{\Box C}.$$
By I.H., we have 
$$\forall x A_0 \dashv \vdash_{\bf IQC} \bigwedge_{C} A^{\Box C}.$$
$\Box$

In the above proof, we do not make use of stability. This means that our proof is simpler than Flagg and Friedman's one in \cite{ff}. Now we can present the following proof.

\it Proof of 2.2\rm . Suppose that $\bf EQC$ $\vdash A^{\Box_{\mathbf{RS}}}$ holds. Let $\Gamma$ be the set of all subformulas of $A$. By Theorems 1.1 and 1.2, we have 
$$\vdash_{\bf IQC} \bigwedge_{C \in \Gamma} (A^{\Box_{\mathbf{RS}}})_\Gamma^{(C)}.$$
Then apply Lemma 2.2 to obtain $\vdash_{\bf IQC} A$. $\Box$

\begin{theorem}  \rm (The soundness and faithfulness for $(\cdot)^{\Box_{\mathbf{RS}}}$) \it For each formula $A$ of $\bf IQC$ $(\bf HA)$, 
$$\vdash_{\bf IQC (HA)} A \enspace \Longleftrightarrow \enspace \enspace \vdash_{\bf EQC (EA)} A^{\Box_{\mathbf{RS}}}.$$
\end{theorem} 

We can also make use of the modified Rasiowa and Sikorski's translation to prove the following conservativeness shown as Corollary 2.1.

\begin{corollary} $\bf EA$ $($$\bf EQC$$)$ is a conservative extension of $\bf HA$$($$\bf IQC$$)$.

\end{corollary} 

Proof. Immediate from Theorem 2.2. $\Box$

\section{Applications of the modified Rasiowa-Sikorski translation}

We recall the following properties.

(1) $\bf EA$ has \it the modal disjunction property \rm $\mathsf{MDP}$ if whenever $\bf EA$ $\vdash \Box A \vee \Box B$, either $\bf EA$ $\vdash \Box A$ or $\bf EA$ $\vdash \Box B$.

(2) $\bf EA$ has \it the modal numerical existence property \rm $\mathsf{MNEP}$ if whenever $\bf EA$ $\vdash \exists x \Box A$, there is some natural number $n$ such $\bf EA$ $\vdash \Box A(\mathbf{n})$ with the numeral $\mathbf{n}$ corresponding to $n$.

(3) $\bf HA$ has \it the disjunction property \rm $\mathsf{DP}$ if whenever $\bf HA$ $\vdash A \vee B$, either $\bf HA$ $\vdash A$ or $\bf HA$ $\vdash B$.

(4) $\bf HA$ has \it the numerical existence property \rm $\mathsf{NEP}$ if whenever $\bf HA$ $\vdash \exists x A$, there is some natural number $n$ such $\bf HA$ $\vdash A(\mathbf{n})$ with the numeral $\mathbf{n}$ corresponding to $n$.

It is well-known that:

\begin{theorem} \rm (Shapiro \cite{shapiro}) \it $\bf EA$ has $\mathsf{MDP}$ and $\mathsf{MNEP}$.
\end{theorem}

For Theorem 3.1, Shapiro \cite{shapiro} proved the version without $\Box$ in conclusion as follows:

$\mathsf{EDP}$: $\bf EA$ $\vdash \Box A \vee \Box B$ $\Longrightarrow$  ($\bf EA$ $\vdash A$ or $\bf EA$ $\vdash B$).

$\mathsf{ENEP}$:  $\bf EA$ $\vdash \exists x \Box A$ $\Longrightarrow$ (there is some natural number $n$ such $\bf EA$ $\vdash A(\mathbf{n})$ with the numeral $\mathbf{n}$ corresponding to $n$).
\smallskip

\noindent Note that $\bf EA$ has the rule of necessitation, of course.
\smallskip

\begin{theorem} \rm (Kleene \cite{kleene45, kleene1}) \it $\bf HA$ has $\mathsf{DP}$ and $\mathsf{NEP}$.
\end{theorem} 

In this section, we shall give an alternative proof of Theorem 3.2. We first observe:

\begin{lemma} $ $

$(1)$ If $\bf EA$ has $\mathsf{MDP}$, then $\bf HA$ has $\mathsf{DP}$.

$(2)$ If $\bf EA$ has $\mathsf{MNEP}$, then $\bf HA$ has $\mathsf{NEP}$.
\end{lemma} 

Proof of (1). We see that

\begin{tabbing}
\hspace*{.4in} \= \hspace{.3in} \= \kill

\>  \> $\bf HA$ $\vdash A \vee B$  \\

\> $\Rightarrow$ \> $\bf EA$ $\vdash (A \vee B)^{\Box_{\mathbf{RS}}}$ \enspace \enspace (Theorem 2.3) \\

\> $\Rightarrow$ \> $\bf EA$ $\vdash \Box A^{\Box_{\mathbf{RS}}} \vee \Box B^{\Box_{\mathbf{RS}}}$ \enspace \enspace (Definition 2.1.(iii)) \\

\> $\Rightarrow$ \> $\bf EA$ $\vdash \Box A^{\Box_{\mathbf{RS}}}$ or $\bf EA$ $\vdash \Box B^{\Box_{\mathbf{RS}}}$ \enspace \enspace ($\mathsf{MDP}$) \\

\> $\Rightarrow$ \> $\bf EA$ $\vdash A^{\Box_{\mathbf{RS}}}$ or $\bf EA$ $\vdash B^{\Box_{\mathbf{RS}}}$ \enspace \enspace ($\mathbf{S4}$) \\

\> $\Rightarrow$ \> $\bf HA$ $\vdash A$ or $\bf HA$ $\vdash B$ \enspace \enspace (Theorem 2.3)  $\Box$

\end{tabbing}

Proof of (2). We see that

\begin{tabbing}
\hspace*{.4in} \= \hspace{.3in} \= \kill

\>  \> $\bf HA$ $\vdash \exists x A(x)$  \\

\> $\Rightarrow$ \> $\bf EA$ $\vdash (\exists x A(x))^{\Box_{\mathbf{RS}}}$ \enspace \enspace (Theorem 2.3) \\

\> $\Rightarrow$ \> $\bf EA$ $\vdash \exists x \Box A(x)^{\Box_{\mathbf{RS}}}$ \enspace \enspace (Definition 2.1.(vi)) \\

\> $\Rightarrow$ \> $\bf EA$ $\vdash \Box A(\mathbf{n})^{\Box_{\mathbf{RS}}}$ \enspace for some $n$ \enspace \enspace ($\mathsf{MNEP}$) \\

\> $\Rightarrow$ \> $\bf EA$ $\vdash A(\mathbf{n})^{\Box_{\mathbf{RS}}}$ \enspace for some $n$ \enspace \enspace ($\mathbf{S4}$) \\

\> $\Rightarrow$ \> $\bf HA$ $\vdash A(\mathbf{n})$ \enspace for some $n$ \enspace \enspace (Theorem 2.3)  $\Box$

\end{tabbing}

The proof of Lemma 3.1 was my original motivation of introducing $(\cdot)^{\Box_{\mathbf{RS}}}$.

The intended alternative proof of Theorem 3.2 is now clear. That is, on the basis of Theorem 3.1 and 
Lemma 3.1, we immediately have Theorem 3.2. 

My suggestion is that the modified Rasiowa and Sikorski's translation would be able to be extened to type theory and set theory. If so, analogous reductions as done in Lemma 3.1 would be possible for meta-theorems of the theories. 

Lastly, we shall give a problem:

\medskip

\noindent \bf Open problem\rm . Is it possible to have the converse reductions of Lemma 3.1 ?

\section{Epistemic Markov's rule}

Unless we specially mention, we shall always follow all formulations and conventions as the above for Heyting arithmetic $\bf HA$, epistemic arithmetic $\bf EA$ (i.e. Peano arithmetic $\bf PA$ based on modal calculus $\bf S4$), intuitionistic first-order calculus with equality $\bf IQC$ and epistemic first-order calculus (i.e. classical first-order calculus based on modal calculus $\bf S4$) with equality $\bf EQC$. By $\bf CQC$, we also denote classical first-order calculus with equality. 

Let $(\cdot)^{(E)}_\Gamma$ be Flagg and Friedman's translation from $\bf EQC$ ($\bf EA$) to $\bf IQC$ ($\bf HA$) (with respect to $\Gamma$ and $E$). ($(\cdot)^{(E)}_\Gamma$ is defined in Definition 1.3.)

The following lemma will make our arguments much simpler.

\begin{lemma}
Let $A$ be a formula of $\bf EQC$. For any finite set $\Gamma$ of formulas of $\bf IQC$ and any $E \in \Gamma$, we have $\vdash_{\bf IQC} \neg_E \neg_E A^{(E)}_\Gamma \leftrightarrow A^{(E)}_\Gamma$.
\end{lemma}

Proof. Easy by induction on the complexity of $A$. $\Box$

Lemma 4.1 is a generalization of the comment in \cite[p. 56]{ff} to that in predicate calculi.

We shall define a translation $(\cdot)^F$ as a simplification of $(\cdot)^{\bot}$ (defined in Definition 1.3) as follows. (Recall that $(\cdot)^{\bot}$ is a special instance of Flagg and Friedman's translation.)

\begin{definition} \rm For any formula $A$ of $\bf EQC$ ($\bf EA$), we inductively define a formula $A^F$ of $\bf IQC$ ($\bf HA$) as follows: 

\begin{tabbing}

\hspace*{.4in} \= \hspace{.4in} \= \hspace{.8in} \= \hspace{.4in} \=  \kill

\> (i) \> $A^F$ \> $=$ \> $\neg^2 A$ for any atomic $A$. \\

\> (ii) \> $(A \vee B)^F$ \> $=$ \> $\neg ^2 (A^F \vee B^F)$. \\

\> (iii) \> $(A \# B)^F$ \> $=$ \> $A^F \# B^F$, $\enspace$ $\# \in \{\wedge, \supset\}$. \\

\> (iv) \> $(\Box A)^F$ \> $=$ \> $A^F$. \\

\> (v) \> $(\exists x A)^F$ \> $=$ \> $\neg^2 \exists x A^F$. \\

\> (vi) \> $(\forall x A)^F$ \> $=$ \> $\forall x A^F$.

\end{tabbing}

\end{definition} 

\begin{proposition}
For any formula $A$ of $\bf EQC$, we have $\vdash_{\bf IQC}A^{\bot} \leftrightarrow A^F$.
\end{proposition}

Proof. In view of Lemma 4.1, it is easy by induction on the complexity of $A$.  $\Box$

\begin{corollary}
For any formula $A$ of $\bf EQC$ $(\bf EA)$, 
$$\vdash_{\bf EQC (EA)} A \enspace \Rightarrow \enspace \vdash_{\bf IQC (HA)} A^F.$$
\end{corollary}

Proof. Obvious from Theorem 1.3 and Proposition 4.1. $\Box$

Let us first recall some elementary facts for $\bf IQC$, which are easily verified.

\begin{proposition}
Let $A$ and $B$ and $E$ be formulas of $\bf IQC$. Then,

\medskip 

$(1)$ \enspace $\vdash_{\bf IQC} A \leftrightarrow \neg_A^2 A$.

\medskip 

$(2)$ \enspace $\vdash_{\bf IQC} \forall x \neg_E^3 A \leftrightarrow \forall x \neg_E A$.

\medskip 

$(3)$ \enspace $\vdash_{\bf IQC} \forall x \neg_E^2 \exists y \neg_E^2 A \leftrightarrow \forall x \neg_E^2 \exists yA$.
\end{proposition}

The following is one of our main theorems.

\begin{theorem} \rm (announced in Inou\'{e} \cite{ino90b}) \it Markov's rule $\mathsf{MR}$ in $\bf HA$ implies epistemic Markov's rule $\mathsf{EMR}$ in $\bf EA$.
\end{theorem}

Proof. We see that

\medskip 

$\vdash_{\bf EA} \Box \forall x \Diamond \exists yR(x, y)$ $\enspace$ (Assumption)

\medskip 

\noindent $\Rightarrow$ $\vdash_{\bf HA} (\Box \forall x \Diamond \exists yR(x, y))^F$ $\enspace$ (Corollary 4.1)

\medskip 

\noindent $\Rightarrow$ $\vdash_{\bf HA} \forall x[\{\neg^2 \exists y \neg^2 R(x, y) \supset \neg^2 \bot\} \supset \neg^2 \bot]$ $\enspace$ (Definition 4.1)

\medskip 

\noindent $\Rightarrow$ $\vdash_{\bf HA} \forall x \neg^4 \exists y \neg^2 R(x, y)$ $\enspace$ (Proposition 4.2.(1))

\medskip 

\noindent $\Rightarrow$ $\vdash_{\bf HA} \forall x \neg^2 \exists y \neg^2 R(x, y)$ $\enspace$ (Proposition 4.2.(2))

\medskip 
	 
\noindent $\Rightarrow$ $\vdash_{\bf HA} \forall x \neg^2 \exists yR(x, y)$ $\enspace$ (Proposition 4.2.(3))

\medskip 
	 
\noindent $\Rightarrow$ $\vdash_{\bf HA} \forall x \exists yR(x, y)$ $\enspace$ (Markov's rule $\mathsf{MR}$)

\medskip 
	 
\noindent $\Rightarrow$ $\vdash_{\bf EA} \forall x \exists yR(x, y)$ $\enspace$ ($\bf EA$ is stronger than $\bf HA$.)

\medskip 
	 
\noindent $\Rightarrow$ $\vdash_{\bf EA} \Box \forall x \exists yR(x, y)$ $\enspace$ (The rule of necessitation). $\Box$

\medskip 

What does Theorem 4.1 mean?  It would be a realization of the constructive content which $\mathsf{EMR}$ has.

Let $(\cdot)^{\Box}$ be G\"{o}del translation (Definition 1.1).

For the last two implications in the above proof of Theorem 4.1, we can apply $(\cdot)^{\Box}$ to $\vdash_{\bf HA} \forall x \exists yR(x, y)$ to directly obtain $\vdash_{\bf EA} \Box \forall x \exists yR(x, y)$.

We need some preparations for the second proof of Theorem 4.1.

\begin{definition} \rm For any formula $A$ of $\bf EQC$ $(\bf EA)$, we define a formula $A^{d\Box}$ of $\bf CQC$ $(\bf PA)$ as follows: 

\begin{tabbing}

\hspace*{.4in} \= \hspace{.4in} \= \hspace{.8in} \= \hspace{.4in} \=  \kill

\> (i) \> $A^{d\Box}$ \> $=$ \> $A$ for any atomic $A$. \\

\> (ii) \> $(A \#  B)^{d\Box} $ \> $=$ \> $A^{d\Box} \#  B^{d\Box}$, $\enspace$ $\# \in \{\vee, \wedge, \supset\}$. \\

\> (iii) \> $(QxA)^{d\Box}$ \> $=$ \> $QxA^{d\Box}$, $Q \in \{\exists,  \forall\}$.

\end{tabbing}

\end{definition}

\begin{proposition}
For any formula $A$ of $\bf EQC$ $(\bf EA)$, 
$$\vdash_{\bf EQC (EA)} A \enspace \Rightarrow\enspace \vdash_{\bf CQC (PA)} A^{d\Box}.$$
\end{proposition}

Proof. Easy by induction on derivation. $\Box$

It is obvious that the translation $(\cdot)^{d\Box}$, i.e. the $\Box$-deleting operator, is not faithful. 

\bigskip

\it The second Proof of Theorem 4.1\rm . We see that

\medskip  

$\vdash_{\bf EA} \Box \forall x \neg \Box \neg \exists yR(x, y)$ $\enspace$ (Assumption)

\medskip 

\noindent $\Rightarrow$ $\vdash_{\bf PA} (\Box \forall x \neg \Box \neg \exists yR(x, y))^{d\Box}$ $\enspace$ (Proposition 4.3)

\medskip 

\noindent $\Rightarrow$ $\vdash_{\bf PA} \forall x \neg^2 \exists yR(x, y)$ $\enspace$ (Definition 4.2)

\medskip 

\noindent $\Rightarrow$ $\vdash_{\bf EA} \forall x \neg^2 \exists yR(x, y)$ $\enspace$ ($\bf EA$ is stronger than $\bf HA$.)

\medskip 

\noindent $\Rightarrow$ $\vdash_{\bf HA} \forall x \neg^2 \exists yR(x, y)$ $\enspace$ (Corollary 2.1).

\medskip 

\noindent Then, we can continue the proof as in the first one. (N.B. In this proof, we may assume, as the crucial point of this proof, that the languages of $\bf HA$, $\bf PA$ and $\bf EA$ contain a common symbol for every primitive recursive function.) $\Box$

\begin{theorem} The following $\mathrm{(1)}$ and $\mathrm{(2)}$ are equivalent.

$\mathrm{(1)}$ Epistemic Markov's rule $\mathsf{EMR}$ holds in $\bf EA$.

$\mathrm{(2)}$ Markov's rule $\mathsf{MR}$ holds in $\bf HA$.
\end{theorem}

Proof. (1) $\Rightarrow$ (2) comes from both Lemma 3.2 and Theorem 3.1 in 
Flagg and Friedman \cite[p. 59]{ff}. 
 And (2) $\Rightarrow$ (1) follows from Theorem 4.1 of this paper. $\Box$

Flagg and Friedman \cite{ff} also directly showed that Epistemic Markov's rule $\mathsf{EMR}$ 
holds in $\bf EA$ as Lemma 3.2 in \cite[p. 59]{ff}. 

\section{Some relations among translations}

In this section, we shall show some relation among the above introduced translations and well-known other ones. 

Let $(\cdot)^\circ$ be G\"{o}del and Gentzen-Bernays's negative translation. We shall recall the translation as follows.

\begin{definition} \rm 
Let $X$ be $\bf CQC$ ($\bf CPC$) and $Y$ $\bf IQC $ ($\bf IPC$). Then we define a translation $(\cdot)^\circ$ from $X$ to $Y$ (called \it G\"{o}del translation\rm ) as follows: for any formula $A$ of $F_X$,

\medskip 

(1) $A^\circ$ = $\neg \neg A$ for $A$ atomic,

\medskip 

(2) $(\neg A)^\circ$ = $\neg A^\circ$,

\medskip 

(3) $(A \vee B)^\circ$ = $\neg (\neg A^\circ \wedge \neg B^\circ$), 

\medskip

(4) $(A \# B)^\circ$ = $A^\circ \# B^\circ$, $\# \in \{\wedge, \supset\}$,

\medskip 

(5) $(\forall x A)^\circ$ = $\forall x A^\circ$,

\medskip 

(6) $(\exists x A)^\circ$ = $\neg \forall x \neg A^\circ$.

\end{definition}

The following theorem is very well-known.

\begin{theorem} \rm (G\"{o}del-Gentzen-Bernays) (\cite{gedel}) \it G\"{o}del translation $(\cdot)^\circ$ is an embedding of $\bf CQC$ $(\bf CPC)$ in $\bf IQC$ $(\bf IPC)$.

\end{theorem} 

Let $(\cdot)^F$ be the translation defined by Definition 4.1. 

\begin{proposition} 
For any formula $A$ of $\bf EQC$, we have 
$$\vdash_{\bf IQC} A^F \leftrightarrow (A^{d\Box})^\circ.$$
\end{proposition}

Proof. Easy by induction on the complexity of $A$.  Recall the following elementary facts: $\vdash_{\bf IQC} \neg \neg (A \vee B)$ $\leftrightarrow$ $\neg (\neg A \wedge \neg B)$ and $\vdash_{\bf IQC} \neg \neg \exists xA$ $\leftrightarrow$ $\neg \forall x \neg A$. $\Box$

From Proposition 5.1, we can understand that Flagg and Friedman's translation contains G\"{o}del and Gentzen-Bernays's negative translation. Also from the same lemma, $(\cdot)^F$ and $((\cdot)^{d\Box})^\circ$ can be, loosely speaking, regarded as the same translations. For any formula $A$ of $\bf CPC$, $A^\circ$ is an invariant under $(\cdot)^\circ$ with respect to the provability in $\bf IQC$ (i.e. $\vdash_{\bf IQC} (A^\circ)^\circ \leftrightarrow A^\circ$). Hence, so is $A^F$ under $(\cdot)^F$ for any formula $A$ of $\bf CPC$, because of Proposition 5.1. Here are another interesting lemma of which the proof is easily proved by induction on the complexity of a given formula.

\begin{lemma} For any formula $A$ of $\bf EQC$, we have:

\medskip 

$(1)$ \enspace If $A$ has no $\forall$, then $\vdash_{\bf IQC} A^F \leftrightarrow \neg \neg A^{d\Box}$.

\medskip 

$(2)$ \enspace $\vdash_{\bf IQC} A^F \leftrightarrow \neg \neg A^F$. 
\end{lemma}

However, we may also make use of Lemma 4.1 and Propositions 4.2 to prove that (2).

\section{Flagg and Friedman's translation and Glivenko's theorem}

Let $\bf CPC$ and $\bf IPC$ be classical and intuitionistic propositional logic, respectively. Glivenko's theorem to be considered in this section is the following and we shall prove it using Flagg and Friedman's translation.

\begin{theorem} \rm (Glivenko \cite{glivenko}) \it For any formula $A$ of $\bf CPC$, we have 
$$\vdash_{\bf CPC} A \enspace \Longleftrightarrow \enspace \vdash_{\bf IPC} \neg \neg A.$$
\end{theorem}

Proof. One direction of the equivalence is trivial. For the other implication, we see that for any formula $A$ of $\bf CPC$, 

\medskip 

$\vdash_{\bf CPC} A$ $\enspace$ (Assumption)

\medskip 

\noindent $\Rightarrow$ $\vdash_{\bf S4} A$ $\enspace$ ($\bf CPC$ is a subsystem of $\bf S4$.)

\medskip 

\noindent $\Rightarrow$ $\vdash_{\bf IPC} A^F$ $\enspace$ (Corollary 4.1)

\medskip 

\noindent $\Rightarrow$ $\vdash_{\bf IPC} \neg \neg A^{d\Box}$ $\enspace$ (Lemma 5.1.(1))

\medskip 

\noindent $\Rightarrow$ $\vdash_{\bf IPC} \neg \neg A$ $\enspace$ ($A$ does not contain $\Box$ !!). \enspace \enspace $\Box$

\section{Horsten's modal-epistemic arithmetic and Markov's rules}

In Shapiro's epistemic arithmetic $\bf EA$, the operator $\Box$ (which is Shapiro's $K$, properly speaking) was introduced as provability (or knowability): that is, $\Box A$ intuitively means that $A$ is provable. But we can further analyse $\Box A$, for example as ``it is possible that we have a proof of $A$". So the provability-operator $\Box$ can be divided to two parts, namely the possibility-operator $\Diamond$ and the proof-operator $P$.  In this context, we can regard $\Box A$ as $\Diamond PA$.  This is the motivation of Leon Horsten to introduce his modal-epistemic arithmetic $\bf MEA$ (see Horsten \cite{horsten1}). In the formulation of $\bf MEA$, $\Box$ and $P$ are taken as primitive operators instead of only $\Box$ in the sense of Shapiro. The possibility-operator  is defined as usual. The Hilbert-style formulation of $\bf MEA$ is the following together with that of Peano arithmetic.

\bigskip

\noindent \sf Modal axioms and rule\rm :
 
\begin{tabbing}
\hspace*{.05in} \= \hspace{.4in} \= \hspace{.5in} \= \kill

\> (i) \> M*1. \> $\vdash_{\bf MEA} \Box A \supset A$. \\

\> (ii) \> M*2. \> $\vdash_{\bf MEA} \Box A \wedge \Box (A \supset B). \supset \Box B$. \\

\> (iii) \> M*3. \> $\vdash_{\bf MEA} A$ $\Rightarrow$ $\vdash_{\bf MEA} \Box A$. \\

\> (iv) \> M*4. \> $\vdash_{\bf MEA} \Diamond A \supset \Box \Diamond A$. \\

\> (v) \> M*5. \> $\vdash_{\bf MEA} \Diamond A \supset A$, where $A$ does not contain $\Box$ and $P$. \\

\> (vi) \> M*6. \> $\vdash_{\bf MEA} \Diamond \exists xA \supset \exists \Diamond A$.

\end{tabbing}

\noindent \sf Epistemic axioms\rm :

\begin{tabbing}
\hspace*{.05in} \= \hspace{.4in} \= \hspace{.5in} \= \kill 

\> (i) \> E*1. \> $\vdash_{\bf MEA} PA \supset A$. \\

\> (ii) \> E*2. \> $\vdash_{\bf MEA} PA \supset PPA$.

\end{tabbing}

\noindent \sf Modal-epistemic axiom and rule\rm : 

\begin{tabbing}
\hspace*{.05in} \= \hspace{.4in} \= \hspace{.5in} \= \kill

\> (i) \> ME*1. \> $\vdash_{\bf MEA} A$ $\Rightarrow$ $\vdash_{\bf MEA} \Diamond PA$. \\

\> (ii) \> ME*2. \> $\vdash_{\bf MEA} \Diamond PA \wedge \Diamond P(A \supset B). \supset \Diamond PB$.

\end{tabbing}

It is not difficult to see that $\bf MEA$ is a conservative extension of Heyting arithmetic $\bf HA$, if we interpret Shapiro's $\Box$ in $\bf EA$ as $P$ (for the details, see \cite[pp. 86--92]{horsten1}): in other words, if we slightly modify the definition of Flagg and Friedman's translation as follows: for any formula of $A$ of $\bf MEA$, (i) $(\Box A)^{(E)}_\Gamma$ = $A^{(E)}_\Gamma$; (ii) $(PA)^{(E)}_\Gamma$ = $\bigwedge_{C \in \Gamma} A^{\Box (C)}_\Gamma$, besides the same definition for the rest of connectives as that of Flagg and Friedman (as a result, we have $(\Diamond PA)^{(E)}_\Gamma$ = $\neg_E \neg_E \bigwedge_{C \in \Gamma} A^{\Box (C)}_\Gamma$, exactly speaking, $\vdash_{\bf HA}$ $(\Diamond PA)^{(E)}_\Gamma$ $\leftrightarrow$ $\neg_E \neg_E \bigwedge_{C \in \Gamma} A^{\Box (C)}_\Gamma$, we can make use of Flagg and Friedman's argument for the faithfulness.

Here we wish to propose several versions of Markov's rules for primitive recursive predicates in $\bf MEA$ as follows:

\medskip 

$\mathsf{HMR_{Nt}}$: $\vdash_{\bf MEA} \Diamond P \forall x \neg \Diamond P \neg \exists y R(x, y)$ $\Rightarrow$ $\vdash_{\bf MEA} \Diamond P \forall x \exists y R(x, y)$.

\medskip 

$\mathsf{HMR_{W1}}$: $\vdash_{\bf MEA} \Diamond \forall x \Diamond P \exists y R(x, y)$ $\Rightarrow$ $\vdash_{\bf MEA} \Diamond P \forall x \exists y R(x, y)$.

\medskip 

$\mathsf{HMR_{W2}}$: $\vdash_{\bf MEA} \Box \forall x \neg \Diamond P \neg \exists y R(x, y)$ $\Rightarrow$ $\vdash_{\bf MEA} \Box \forall x \exists y R(x, y)$,

\medskip 

$\mathsf{HMR_{W3}}$: $\vdash_{\bf MEA} \Box \forall x \Diamond \exists y R(x, y)$ $\Rightarrow$ $\vdash_{\bf MEA} \Box \forall x \exists y R(x, y)$.

\medskip 

$\mathsf{HMR_{W4}}$: $\vdash_{\bf MEA} P \forall x \neg \Diamond P \neg \exists y R(x, y)$ $\Rightarrow$ $\vdash_{\bf MEA} P \forall x \exists y R(x, y)$.

\medskip 

$\mathsf{HMR_{W5}}$: $\vdash_{\bf MEA} P \forall x \Diamond \exists y R(x, y)$ $\Rightarrow$ $\vdash_{\bf MEA} P \forall x \exists y R(x, y)$.

\medskip 

\noindent Hereby $R(x, y)$ is $f (x, y) = 0$ for some binary primitive recursive function symbol $f$. It is easy to see that 

\bigskip

(i) \enspace $\mathsf{HMR_{Nt}}$ $\Leftrightarrow$ $\mathsf{HMR_{W2}}$.

\medskip

(ii) \enspace $\mathsf{HMR_{W1}}$ $\Rightarrow$ $\mathsf{HMR_{W3}}$.

\begin{theorem}
$\bf MEA$ is closed under any of $\mathsf{HMR_{Nt}}$ , $\mathsf{HMR_{W1}}$, $\mathsf{HMR_{W2}}$ and $\mathsf{HMR_{W3}}$.
\end{theorem}

Proof. Similar to that of Lemma 3.2 in \cite{ff}. In addition, it is also easy to have a direct proof in $\bf MEA$ of $\mathsf{HMR_{W1}}$ and $\mathsf{HMR_{W3}}$, respectively. $\Box$

As the proof of Theorem 3.3 in \cite{ff}, using $\mathsf{HMR_{Nt}}$, we have an alternative proof of:

\begin{theorem}
$\bf HA$ is closed under $\mathsf{MR}$.
\end{theorem}

Theorem 7.2 together with its proof tells us that $\mathsf{HMR_{Nt}}$ is the most natural candidate of Markov's rule for $\bf MEA$. 

As shown in Theorem 4.2, the proof in $\bf MEA$ of Theorem 7.2 can be partly reduced to that in $\bf HA$ as follows:

\begin{theorem}
$\mathsf{MR}$ implies $\mathsf{HMR_{Nt}}$ , $\mathsf{HMR_{W1}}$, $\mathsf{HMR_{W2}}$ and $\mathsf{HMR_{W3}}$.
\end{theorem}

Proof. Similar to that of Theorem 4.2. $\Box$

In addition, epistemic Markov's principle has been to some extent considered in Horsten \cite{horsten1, horsten93b}. On the contrary, to the best of my knowledge, there has been so far no study of general epistemic Markov's Rule. The investigation on epistemic mathematics and its metamathematics is still in the age of pioneers.

Lastly, we wish to pose the following question: is it possible to give a Hilbert-style formulation of $\bf MEA$ without modal-epistemic axioms and rules so that its modal notion and its epistemic one completely are separated with each other?

\section{A proof of Fundamental Conjecture $\mathsf{FC}$ for Friedman's theorem}

In this section, we shall give a partial solution to Fundamental Conjecture $\mathsf{FC}$ (i.e. $\mathsf{FS}$ $\Rightarrow$ $\mathsf{F}$) which was already mentioned in Introduction \S 1.

Let us recall the specification of notations to be used in this section (see \S 3.5 for the details): 

\medskip

\noindent $\mathsf{DP}$ = The disjunction property in $\bf HA$ and its non-modal extensions.

\medskip

\noindent $\mathsf{NEP}$ = The numerical existence property in $\bf HA$ and its non-modal extensions.

\medskip

\noindent $\mathsf{MDP}$ = The modal disjunction property in a modal system.

\medskip

\noindent $\mathsf{MNEP}$ = The modal numerical existence property in a modal system.

\medskip

\noindent $\mathsf{EDP}$ = The modal disjunction property in $\bf EA$.

\medskip

\noindent $\mathsf{ENEP}$ = The modal numerical existence property in $\bf EA$.

\medskip

First we shall recall the following (Lemma 3.1):

\medskip

$(1)$ \enspace $\mathsf{EDP}$ $\Rightarrow$ ($\mathsf{DP}$ in $\bf HA$).

\medskip

$(2)$ \enspace $\mathsf{ENEP}$ $\Rightarrow$ ($\mathsf{NEP}$ in $\bf HA$).

\medskip

We also note:

\begin{theorem} In $\bf HA$, we have
$$\mathsf{NEP} \enspace \Rightarrow \enspace \mathsf{DP}.$$
\end{theorem}

Proof. It follows from
$$\vdash_{\bf HA} \enspace A \vee B. \leftrightarrow \exists x(x = 0 \supset A. \wedge . x \neq 0 \supset B).$$
$\Box$

In the description of $\mathsf{FC}$, by $\mathsf{FS}$ (= a modal analogue of Friedman's theorem, i.e. Theorem in \cite{fs}), we mean the following.

\begin{theorem} \rm 
(The corrected version of Friedman and Sheard \cite{fs}). \it Suppose that $S$ is a recursively enumerable theory including $\bf PA$ $($N.B. the language of S is a first-order language appropriate for arithmetic, argumented by the modal operator $\Box$$)$ such that

$(1)$ \enspace $\vdash_S \Box A \supset \Box \Box A$ for any formula $A$.

$(2)$ \enspace $\vdash_S A \supset \Box A$ for any $\Delta_0$ formula $A$.

$(3)$ \enspace $\vdash_S \Box (\Box A \supset A)$ for any $\Delta_0$ formula $A$.

$(4)$ \enspace If $\vdash_S \Box A$, then  $\vdash_S A$ for any formula $A$, 

$(5)$ \enspace If $A_1, \dots , A_n \vdash_S B$, then  $\vdash_S (\Box A_1 \wedge \cdots \wedge \Box A_n) \supset  \Box B$ for arbitrary formulas $A_1, \dots , A_n, B$ $(n \geq 0)$ of first-order logic, 

\noindent where $\Delta_0$ formulas are the usual formulas of the arithmetical hierachy with only bounded quantifiers and with no occurrences of $\Box$. Then, $S$ satisfies $\mathsf{MDP}$ if and only if $S$ satisfies $\mathsf{MNEP}$.
\end{theorem}

 Friedman and Sheard \cite{fs} noted the following theorem.
 
 \begin{theorem} \rm 
(The corrected version of Friedman and Sheard \cite{fs}).  \it Any recursively enumerable extension of $\bf EA$ which has $\mathsf{MDP}$ also has $\mathsf{MNEP}$.
\end{theorem}

By $\mathsf{FS_{\bf EA}}$ we denote Theorem 8.3.

In the description of $\mathsf{FC}$, by $\mathsf{F}$ (= Friedman's theorem, i.e. Theorem 1 in \cite{friedman}), we mean the following one.

\begin{theorem} \rm 
(Friedman \cite{friedman}; cf. Leivant \cite{leivant}). \it Any recursively enumerable extension of $\bf HA$ which has $\mathsf{DP}$ also has $\mathsf{NEP}$.
\end{theorem}

The main theorem of this section is the affirmative answer of our \it Fundamental Conjecture \rm (for short, $\mathsf{FC}$):
 $$\mathsf{FS} \enspace \Longrightarrow \enspace  \mathsf{F}.$$
 
\noindent $\mathsf{FC}$ was conjectured in Inou\'{e} \cite{ino90a}.

\begin{theorem} The Fundamental Conjecture $\mathsf{FC}$ holds.
\end{theorem}

Proof. Suppose $\mathsf{FS}$ holds. Let a recursively enumerable extension of $\bf HA$ be $T$. 
Take $S$ as  a recursively enumerable theory including $\bf PA$ (The language of $S$ is a first-order language appropriate for arithmetic, argumented by the modal operator $\Box$) such that

$(1)$ \enspace $\vdash_S \Box A \supset \Box \Box A$ for any formula $A$.

$(2)$ \enspace $\vdash_S A \supset \Box A$ for any $\Delta_0$ formula $A$.

$(3)$ \enspace $\vdash_S \Box (\Box A \supset A)$ for any $\Delta_0$ formula $A$.

$(4)$ \enspace If $\vdash_S \Box A$, then  $\vdash_S A$ for any formula $A$, 

$(5)$ \enspace If $A_1, \dots , A_n \vdash_S B$, then  $\vdash_S (\Box A_1 \wedge \cdots \wedge \Box A_n) \supset  \Box B$ for arbitrary formulas $A_1, \dots , A_n, B$ $(n \geq 0)$ of first-order logic.

$(6)$ \enspace $S$ extends $T$. 

Further extend $S$ to $S_{\mathsf{FS}}$ replacing $\bf PA$ by $\bf EA$. 
$S_{\mathsf{FS}}$ is satisfied with the assumption of $\mathsf{FS}$ (See Friedman and Sheard \cite{fs}). So we have, in the language of $S_{\mathsf{FS}}$, the equivalence between $\mathsf{MDP}$ and $\mathsf{MNEP}$. This fact contains the case of formulas translated by modified Rasiowa-Sikorski (or G\"{o}del) translation. Since $\mathsf{FS}$ holds, we have a proof $\mathcal{P(S_{FS}})$ of this translated-formula-case from the proof of $\mathsf{FS}$, adapting the arguments of G\"{o}del number to those of the case, also. Delete all $\Box$'s from $\mathcal{P(S_{FS}})$. Then we obtain the proof the equivalence between $\mathsf{DP}$ and $\mathsf{NEP}$. This completes the proof. $\Box$
\medskip

\it I think that the proof of Theorem 8.5 is an example of a new type of proofs. That is, we make use of the proof itself of a theorem to prove a theorem\rm . This relates to the formalization of mathematics. In order to obtain more rigorous proof of Theorem 8.5, after formalizing $\mathsf{FS}$ in Mizar, then delete the $\Box$'s from the formalized proof and reconstruct the usual proof of $\mathsf{F}$ from it. I shall touch on the idea of this way of proof in mathematics in Inou\'{e} \cite{Inoue2023hum}. For Mizar, refer to Bancerek et al. \cite{ban1, ban2} and Matuszewski et al. \cite{mr}.
\smallskip

We shall give two open problems as follows.
\medskip

\bf Open problem 1. \rm Prove the versions of modal Markov's rule for the modal arithmetics proposed in  Kurahashi and Okuda \cite{kuranishi1} and study the relationships between the modal rule and the non-modal one.
\smallskip

\bf Open problem 2. \rm  In order to obtain more rigorous proof of Theorem 8.5, after formalizing $\mathsf{FS}$ in Mizar, then delete the $\Box$'s from the formalized proof and reconstruct the usual proof of $\mathsf{F}$ from it. The formalization of embedding will be also needed.

\section{Discussion and my Philosophy}

In Inou\'{e} \cite{ino90a}, I wrote my opinion at that time as follows: epistemic systems of a mathematical theory will give us one of sound grounds on which philosophy, mathematical logic and mathematics can be desirably integrated and interchanged. Even though the opinion would sound good for philosophers, most mathematicians in pure and applied mathematics and recent many (but not all) mathematical logicians (= mathematicians) will not be satisfied with it. I believe that most mathematicians expect real mathematical advances in mathematics, when a new theory was introduced. In this sense, studies on epistemic systems will not \it directly \rm contribute mathematics in the usual sense, although a lot of technical results on epistemic systems are of course mathematical ones. Because the character of them is more or less philosophical in such a sense that we can be conscious of the constructiveness lain in mathematics by doing research on epistemic systems. But the study of epistemic systems can well \it \rm contribute mathematics as a subject including recent highly developed mathematical logic, though in an \it indirect \rm way. That is, not all but in such a way as \it changing our recoginition \rm and \it reconsiderations about what was already established\rm . This would be able to become a great impetus to mathematics and mathematicians. 
The same effect can be seen in one of the roles of non-standard analysis (Robinson \cite{robinson1974}, Loeb and Wolff \cite{loebWolff} etc).  However the role of epistemic mathematics is up to the moment very weeker than that of non-standard analysis. We have to find ways of the significance of epistemic mathematics in this direction. Even if we recognize epistemic mathematics as a humanization of mathematics in the sense of Inou\'{e}  \cite{Inoue2023hum}, this direction is still important. Epistemic mathematics should be recognized as a humanization of mathematics. I think that \it recognition determines human intended behavior. Recognition will change the world. \rm That is my Philosophy. 
In this sense, epistemic mathematics is meaningful. 
As another direction of thinking, my recent practices of  philosophy is in Inou\'e \cite{Inoue2022b}, 
which will change our notion of coordinates in physics, and even mathematics and science. The role of epistemic mathematics will be enhanced in Inou\'{e} \cite{Inoue2023modal} in the future. For  \cite{Inoue2022b, Inoue2023hum, Inoue2023modal}, I here only mention the existence of the papers. 
I will in the future develop my thought further, which involves physics and science that contains data 
science (See further  Inou\'{e} \cite{Inoue2023philodata, Inoue2023missing, Inoue2023humphysics}). 

\it Last but not least, I think that the concept of modality and spacial concept will be theoretically closer in the future, so in the sense epistemic mathematics and its derived forms will contribute mathematics more substantially\rm .

\medskip

$$ $$

Takao Inou\'{e}

\medskip

AI Nutrition Project,

Artificial Intelligence Center for Health and Biomedical Research,

National Institutes of Biomedical Innovation, Health and Nutrition 

 (NIBIOHN), Osaka, Japan

\bigskip 

t.inoue@nibiohn.go.jp

(Personal) takaoapple@gmail.com

I prefer my personal email.


\begin{thebibliography}{99}
\bibitem{Inoue2023philodata} \sc T. Inou\'e\rm , {Force fields potentials and images: two core immanent analogies for a novel methodology of data science which uses missing data imputations.} This is a temporal title. Forthcoming.
\bibitem{Inoue2023missing} \sc T. Inou\'e\rm , {Data tables as images, a novel imputation method VOLIM for missing data with a case study of certain real data.} This is a temporal title. Forthcoming.
\bibitem{Inoue2023OnoKleene} \sc T. Inou\'e\rm ,  {A constructive proof of the intuitionistic analogue of Ono-Kleene's theorem and non-C-equivalence.} This is a temporal title. Forthcoming.
\bibitem{Inoue2023humphysics} \sc T. Inou\'e\rm , {On a humanization of physics and science: quantum mechanics and beyond.} This is a temporal title. Forthcoming.
\end{thebibliography}
\end{document}